\documentclass[10pt]{llncs}
\usepackage[]{amsmath}
\usepackage[]{amssymb}

\newcommand{\abs}[1]{{\left|#1\right|}}

\newcommand{\RP}{\mathop\mathsf{RP}}
\newcommand{\PR}{\mathop\mathsf{PR}}
\newcommand{\RPPR}{\mathop\mathsf{RPPR}}
\newcommand{\ANY}{\mathop\mathsf{ANY}}
\newcommand{\ORD}{\mathop\mathsf{ORD}}

\begin{document}

\title{\emph{Addenda/Corrigenda}  \\
Fixed Points and Two-cycles of the Discrete Logarithm}
\author{Joshua Holden}
\institute{Department of Mathematics,
Rose-Hulman Institute of Technology,
Terre Haute, IN, 47803-3999, USA,
\email{holden@rose-hulman.edu}}

\maketitle
\begin{abstract}
    This report consists of additions and corrections to the author's
    paper~\cite{Holden02}, which appeared in the proceedings of the ANTS~V
    conference.  The work described here was presented at the
    conference itself, which took place after the original paper was
    published.  The abstract of the original paper was as follows: We
    explore some questions related to one of Brizolis: does every
    prime $p$ have a pair $(g,h)$ such that $h$ is a fixed point for
    the discrete logarithm with base $g$?  We extend this question to
    ask about not only fixed points but also two-cycles.  Campbell and
    Pomerance have not only answered the fixed point question for
    sufficiently large $p$ but have also rigorously estimated the
    number of such pairs given certain conditions on $g$ and $h$.  We
    attempt to give heuristics for similar estimates given other
    conditions on $g$ and $h$ and also in the case of two-cycles.
    These heuristics are well-supported by the data we have collected,
    and seem suitable for conversion into rigorous estimates in the
    future.
\end{abstract}



The formulas leading up to Conjecture~\ref{conj7} have several typos.
The corrected formulas should read:

$$\sum_{m \mid p-1}\abs{S_{m}}^{2}/\abs{T_{m}} \approx
\sum_{m \mid p-1} \frac{\phi(m)}{m^{2}} \left( \sum_{d \mid (p-1)/m}
\frac{\phi(dm)}{d}\right)^{2}$$

$$N_{(\ref{ha}), h \ANY, a \ANY}(p) \approx (p-1) +
\sum_{m \mid p-1} \frac{\phi(m)}{m^{2}} \left( \sum_{d \mid (p-1)/m}
\frac{\phi(dm)}{d}\right)^{2}$$

(This is also the formula which should appear in
Conjecture~\ref{conj7}(\ref{conj7a}).)

$$\sum_{m \mid p-1} \frac{\phi(m)^{3}}{m^{2}} \left(\sum_{d \mid (p-1)/m}
\frac{\phi(d)}{d} \right)^{2}
= \sum_{m \mid p-1} \left(\prod_{q \mid m}
\frac{\phi(q)^{3}}{q^{2}} \right) \left(\prod_{q \mid (p-1)/m} \left( 1 +
\frac{\phi(q)}{q} \right)^{2} \right)$$
$$= \prod_{q \mid p-1} \left( \frac{\phi(q)^{3}}{q^{2}} + \left( 1 +
\frac{\phi(q)}{q} \right)^{2} \right)
= \prod_{q \mid p-1} \left( q + 1 -
\frac{1}{q} \right)$$

In addition, there are two errors in~(\ref{geneq}); the equation should
    read:
    \begin{multline} \tag{\ref{geneq}}
\sum_{m \mid p-1}\abs{S_{m}}^{2}/\abs{T_{m}} \approx
\prod_{q} \left(
\sum_{\beta=0}^\alpha \phi(q^\beta)
\left[\left(1-\frac{1}{q}\right)(\alpha-\beta) +
\frac{\phi(q^\beta)}{q^\beta} \right]^2\right) \\
\shoveleft{= \prod_q \left(
\left[\left(1-\frac{1}{q}\right)\alpha + 1\right]^2  + \sum_{\beta=1}^\alpha
q^\beta\left(1-\frac{1}{q}\right) \left[
\left(1-\frac{1}{q}\right)(\alpha-\beta+1) \right]^2 \right)}\\
\shoveleft{= \prod_q \left(
\left[\left(1-\frac{1}{q}\right)\alpha + 1\right]^2 \right.} \\
+ \left(1-\frac{1}{q}\right)^3
\left[ (\alpha+1)^2 \frac{q^{\alpha+1}-q}{q-1}
- 2 (\alpha+1) \frac{\alpha q^{\alpha+2} -
(\alpha+1)q^{\alpha+1}+q}{(q-1)^2} \right.\\
\left.\left.
+ \frac{\alpha^2 q^{\alpha+3} - (2\alpha^2 + 2\alpha-1)q^{\alpha+2}
        + (\alpha^2+2\alpha+1)q^{\alpha+1}-q^2-q}{(q-1)^3}
\right]
\right)
\end{multline}

    Conjecture~\ref{conj1}(\ref{conj1c}) is incorrect.
    In~(\ref{fp}) if $h \PR$ then $g \PR$ also, so
    $N_{(\ref{fp}), g \ANY, h \PR}(p)$ is in fact equal to
    $N_{(\ref{fp}), g \PR, h \RPPR}(p)$.


The same observation for~(\ref{tc}) gives $N_{(\ref{tc}), g
\PR, h \RPPR}(p) = N_{(\ref{tc}), g \ANY, h \RPPR}(p)$ in
Conjecture~\ref{conj5} and $N_{(\ref{tc}), g
\PR, h \PR}(p) = N_{(\ref{tc}), g \ANY, h \PR}(p)$ in
Conjecture~\ref{conj6}(\ref{conj6b}).

    In Conjectures~\ref{conj2} and~\ref{conj6}(\ref{conj6a}) it
    should be noted that the observed values in question must be
    exactly equal, by the symmetry of~(\ref{ha}).  Likewise in
    Conjectures~\ref{conj4} and~\ref{conj6}(\ref{conj6b}).

    A complete and corrected list of Theorems and Conjectures
    follows.  The numbering from the original paper has been
    preserved as much as possible.

    \begin{proposition}
    $N_{(\ref{fp}), g \ANY,
h \RP}(p) = \phi(p-1).$
    \end{proposition}

    \begin{theorem}[Zhang, independently by others]
\begin{equation*}
\begin{split}
    N_{(\ref{fp}), g \PR, h \RPPR}(p)& =  N_{(\ref{fp}), g \PR, h
    \RP}(p) \\
    & = N_{(\ref{fp}), g \PR, h \PR}(p) \\
    & = N_{(\ref{fp}), g \ANY, h \RPPR}(p) \\
    & = N_{(\ref{fp}), g \ANY, h \PR}(p) \\
    &\approx
\mbox{$\phi(p-1)^{2}/(p-1)$}
\end{split}
\end{equation*}
\end{theorem}

    \begin{conjecture} \label{conj1}
    \mbox{}
    \begin{enumerate}
\item $N_{(\ref{fp}), g \ANY, h \ANY}(p) \approx
p-1$
\item $N_{(\ref{fp}), g \PR, h \ANY}(p) \approx \phi(p-1)$
\item \label{conj1c}
(See above.)
\item $N_{(\ref{fp}), g \RP, h \bullet}(p)
\approx \phi(p-1)/(p-1) N_{(\ref{fp}), g \ANY, h \bullet}(p)$
\item $N_{(\ref{fp}), g \RPPR, h \bullet}(p)
\approx \phi(p-1)/(p-1) N_{(\ref{fp}), g \PR, h \bullet}(p)$
\end{enumerate}

\end{conjecture}

\begin{conjecture} \label{conj2}
    \begin{equation*}
    \begin{split}
    N_{(\ref{tc}), g \ANY, h \RP}(p) &= N_{(\ref{ha}), h \RP,
a \ANY}(p) \\
&\approx 2 \phi(p-1).
\end{split}
\end{equation*}
\end{conjecture}

\begin{conjecture} \label{conj3}
        \begin{equation*}
    \begin{split}
 N_{(\ref{tc}), h
\RP, g \ORD h}(p)&= N_{(\ref{ha}), h \RP, a \RP}(p) \\
&\approx \phi(p-1) + \mbox{$\phi(p-1)^{2}/(p-1)$}.
\end{split}
\end{equation*}

\end{conjecture}

\begin{conjecture} \label{conj4}
         \begin{equation*}
    \begin{split}
   N_{(\ref{tc}), g \PR, h \RP}(p) &= N_{(\ref{ha}), h
\RP, a \PR}(p) \\
&\approx 2\phi(p-1)^{2}/(p-1).
\end{split}
\end{equation*}

\end{conjecture}

\begin{conjecture} \label{conj5}
         \begin{equation*}
    \begin{split}
    N_{(\ref{tc}), g
\PR, h \RPPR}(p)  &= N_{(\ref{tc}), g \ANY, h \RPPR}(p) \\
&= N_{(\ref{ha}), h \RPPR, a \bullet}(p) \\
&= N_{(\ref{ha}), h \bullet, a \RPPR}(p) \\
&= N_{(\ref{ha}), h \RPPR, a \RPPR}(p) \\
&\approx
\phi(p-1)^{2}/(p-1) + \phi(p-1)^{3}/(p-1)^{2}.
\end{split}
\end{equation*}
\end{conjecture}

\begin{conjecture} \label{conj6} \mbox{}
    \begin{enumerate}
\item \label{conj6a}
         \begin{equation*}
    \begin{split}
N_{(\ref{tc}), h \ANY, g \ORD h}(p) &= N_{(\ref{ha}), h \ANY, a
 \RP}(p) \\
 &\approx 2 \phi(p-1).
\end{split}
\end{equation*}
\item \label{conj6b}
         \begin{equation*}
    \begin{split}
N_{(\ref{tc}), g
\PR, h \PR}(p) &= N_{(\ref{tc}), g \ANY, h \PR}(p) \\
&= N_{(\ref{ha}), h \PR, a \RP}(p) \\
&\approx
2\phi(p-1)^{2}/(p-1).
\end{split}
\end{equation*}
\end{enumerate}
\end{conjecture}

\begin{conjecture} \label{conj7} \mbox{}
    \begin{enumerate}
\item \label{conj7a}
$N_{(\ref{ha}), h \ANY, a \ANY}(p) \approx (p-1) +
\sum_{m \mid p-1} \frac{\phi(m)}{m^{2}} \left( \sum_{d \mid (p-1)/m}
\frac{\phi(dm)}{d}\right)^{2}$.
\item \label{conj7b} If $p-1$ is squarefree then $N_{(\ref{ha}), h \ANY, a \ANY}(p) \approx (p-1) +
\prod_{q \mid p-1} \left( q + 1 -
\frac{1}{q} \right)$, where the product is taken over primes $q$
dividing $p-1$.
\item \label{conj7c} In general, $N_{(\ref{ha}), h \ANY, a \ANY}(p) \approx (p-1)$ plus
the formula given in~(\ref{geneq}).
\item \label{conj7d} $N_{(\ref{ha}), h \PR,
a \ANY}(p) \approx 2 \phi(p-1)$.
\item \label{conj7e}
$N_{(\ref{ha}), h \ANY, a \PR}(p)\approx 2 \phi(p-1)$.
\item \label{conj7f}
 $N_{(\ref{ha}), h \PR, a \PR}(p) \approx \phi(p-1) +
\phi(p-1)^{2}/(p-1)$.

\end{enumerate}
\end{conjecture}
\begin{conjecture} \label{conj8} \mbox{}
    \begin{enumerate}
    \item
$N_{(\ref{tc}), g \PR, h \ANY}(p) \approx 2\phi(p-1)$.
\item
    $N_{(\ref{tc}), g \ANY, h \ANY}(p) \approx
2(p-1)$.
\end{enumerate}
\end{conjecture}

\begin{conjecture} \label{conj9} \mbox{} 
\begin{enumerate}
    \item $N_{(\ref{tc}), g \RP, h
    \bullet}(p) \approx \phi(p-1)/(p-1) N_{(\ref{tc}), g \ANY, h \bullet}(p)$.
\item $N_{(\ref{tc}), g \RPPR, h \bullet}(p)
\approx \phi(p-1)/(p-1) N_{(\ref{tc}), g \PR, h \bullet}(p)$.
\end{enumerate}
\end{conjecture}

These conjectures are summarized in Tables~\ref{fptalktable},
\ref{hatalktable}, and~\ref{tctalktable}, which also contain new data
presented at ANTS~V.  This data was collected on a  Beowulf cluster,
with 19 nodes, each consisting of 2 Pentium III processors running at
1~Ghz.  The programming was done in C, using MPI, OpenMP, and OpenSSL
libraries.  The collection took 68 hours for all values of
$N_{(\bullet), \bullet, \bullet}(p)$, for five primes $p$ starting at 100000.
Table~\ref{relationstable} summarizes the relationships between
solutions to~(\ref{tc}) and solutions to~(\ref{ha}).

\begin{table}[!ht]
    \caption{Solutions to~(\ref{fp})}
    \label{fptalktable}
    $$\begin{array}{|c|c|c|c|c|}
\multicolumn{5}{l}{\text{(a) Predicted formulas for $N_{(\ref{fp})}(p)$}}\\
    \hline
        g \setminus h & \ANY & \PR & \RP & \RPPR  \\
        \hline
    \ANY &\approx \scriptstyle(p-1) &
    \approx\frac{\phi(p-1)^{2}}{(p-1)} &
        = \scriptstyle \phi(p-1) &
    \approx\frac{\phi(p-1)^{2}}{(p-1)} \\

        \hline
    \PR & \approx \scriptstyle \phi(p-1) &
    \approx\frac{\phi(p-1)^{2}}{(p-1)} &
    \approx\frac{\phi(p-1)^{2}}{(p-1)} &
    \approx\frac{\phi(p-1)^{2}}{(p-1)} \\

    \hline
    \RP & \approx \scriptstyle\phi(p-1) &
    \approx\frac{\phi(p-1)^{3}}{(p-1)^2} &
    \approx\frac{\phi(p-1)^{2}}{(p-1)} &
    \approx\frac{\phi(p-1)^{3}}{(p-1)^2} \\

    \hline
    \RPPR & \approx\frac{\phi(p-1)^{2}}{(p-1)} &
    \approx\frac{\phi(p-1)^{3}}{(p-1)^2} &
    \approx\frac{\phi(p-1)^{3}}{(p-1)^2} &
    \approx\frac{\phi(p-1)^{3}}{(p-1)^2} \\
        \hline
\multicolumn{5}{l}{}\\
    \multicolumn{5}{l}{\text{(b) Predicted values for
    $N_{(\ref{fp})}(100057)$}}\\
        \hline
        g \setminus h & \ANY & \PR & \RP & \RPPR  \\
        \hline
\ANY & 100056 & 9139.46 & 30240 & 9139.46 \\
    \hline
    \PR  &30240  &9139.46  &9139.46 & 9139.46 \\
    \hline
    \RP  &30240  &2762.23 & 9139.46 & 2762.23 \\
    \hline
    \RPPR  &9139.46&  2762.23 & 2762.23 & 2762.23  \\
    \hline
\multicolumn{5}{l}{}\\
    \multicolumn{5}{l}{\text{(c) Observed values for
    $N_{(\ref{fp})}(100057)$}}\\

        \hline
        g \setminus h & \ANY & \PR & \RP & \RPPR  \\
        \hline
   \ANY&  98506& 9192 & 30240& 9192\\
   \hline
\PR&   29630& 9192&  9192&  9192\\
\hline
\RP&   29774& 2784&  9037&  2784\\
\hline
\RPPR& 9085&  2784&  2784&  2784\\
\hline
\end{array}$$
\end{table}

\begin{table}[!ht]
    \caption{Solutions to~(\ref{ha})}
    \label{hatalktable}

    $$\begin{array}{|c|c|c|c|c|}

\multicolumn{5}{l}{\text{(a) Predicted formulas for the nontrivial
part of $N_{(\ref{ha})}(p)$}}\\
    \hline
        a \setminus h & \ANY & \PR & \RP & \RPPR  \\
        \hline
    \ANY & \approx
    \sum\frac{\abs{S_{m}}^{2}}{\abs{T_{m}}}  &
    \approx \scriptstyle \phi(p-1) &
    \approx \scriptstyle \phi(p-1) &
    \approx \frac{\phi(p-1)^{3}}{(p-1)^{2}} \\

        \hline
    \PR & \approx \scriptstyle \phi(p-1) &
    \approx \frac{\phi(p-1)^{2}}{(p-1)} &
    \approx \frac{\phi(p-1)^{2}}{(p-1)} &
    \approx \frac{\phi(p-1)^{3}}{(p-1)^{2}} \\

    \hline
    \RP &  \approx \scriptstyle \phi(p-1) &
    \approx \frac{\phi(p-1)^{2}}{(p-1)} &
    \approx \frac{\phi(p-1)^{2}}{(p-1)} &
    \approx \frac{\phi(p-1)^{3}}{(p-1)^{2}} \\

    \hline
    \RPPR & \approx \frac{\phi(p-1)^{3}}{(p-1)^{2}} &
    \approx \frac{\phi(p-1)^{3}}{(p-1)^{2}} &
    \approx \frac{\phi(p-1)^{3}}{(p-1)^{2}} &
    \approx \frac{\phi(p-1)^{3}}{(p-1)^{2}} \\

        \hline
\multicolumn{5}{l}{}\\
\multicolumn{5}{l}{\text{(b) Predicted values for the nontrivial
part of $N_{(\ref{ha})}(100057)$}}\\
        \hline
        a \setminus h & \ANY & \PR & \RP & \RPPR  \\
        \hline
\ANY & 190822.0 & 30240 & 30240 & 2762.225  \\
\hline
\PR & 30240 & 9139.458 & 9139.458 & 2762.225  \\
\hline
\RP & 30240 & 9139.458 & 9139.458 & 2762.225  \\
\hline
\RPPR & 2762.225 & 2762.225 & 2762.225 & 2762.225 \\
    \hline
\multicolumn{5}{l}{}\\
\multicolumn{5}{l}{\text{(c) Observed values for the nontrivial
part of $N_{(\ref{ha})}(100057)$}}\\
        \hline
        a \setminus h & \ANY & \PR & \RP & \RPPR  \\
        \hline
\ANY &  190526 & 30226 & 30291 & 2820 \\
\hline
\PR &   30226 & 9250 & 9231 & 2820 \\
\hline
\RP &   30291 & 9231 & 9086 & 2820\\
\hline
\RPPR & 2820 & 2820 & 2820 & 2820 \\
\hline
\end{array}$$
\end{table}

\begin{table}[!ht]
    \caption{Solutions to~(\ref{tc})}
    \label{tctalktable}

   $$\begin{array}{|c|c|c|c|c|}
\multicolumn{5}{l}{\text{(a) Predicted formulas for the nontrivial
part of $N_{(\ref{tc})}(p)$}}\\
    \hline
        g \setminus h & \ANY & \PR & \RP & \RPPR  \\
        \hline
    \ANY & \approx \scriptstyle (p-1)    &
    \approx \scriptstyle \phi(p-1)  &
    \approx \frac{\phi(p-1)^{2}}{(p-1)} &
    \approx \frac{\phi(p-1)^{3}}{(p-1)^{2}} \\

        \hline
    \PR & \approx \scriptstyle \phi(p-1)   &
    \approx \frac{\phi(p-1)^{2}}{(p-1)} &
    \approx \frac{\phi(p-1)^{2}}{(p-1)} &
    \approx \frac{\phi(p-1)^{3}}{(p-1)^{2}} \\

    \hline
    \RP &  \approx \scriptstyle \phi(p-1) &
    \approx \frac{\phi(p-1)^{2}}{(p-1)} &
    \approx \frac{\phi(p-1)^{3}}{(p-1)^{2}} &
    \approx \frac{\phi(p-1)^{4}}{(p-1)^{3}} \\

    \hline
    \RPPR & \approx \frac{\phi(p-1)^{2}}{(p-1)}     &
    \approx \frac{\phi(p-1)^{3}}{(p-1)^{2}} &
    \approx \frac{\phi(p-1)^{3}}{(p-1)^{2}} &
    \approx \frac{\phi(p-1)^{4}}{(p-1)^{3}} \\

        \hline
\multicolumn{5}{l}{}\\
\multicolumn{5}{l}{\text{(b) Predicted values for the nontrivial
part of $N_{(\ref{tc})}(100057)$}}\\
        \hline
        g \setminus h & \ANY & \PR & \RP & \RPPR  \\
        \hline
\ANY & 100056 & 9139.5 & 30240 & 2762.2  \\
\hline
\PR &30240 & 9139.5 & 9139.5 & 2762.2 \\
\hline
\RP & 30240 & 2762.2 & 9139.5 & 834.8   \\
\hline
\RPPR & 9139.5 & 2762.2 & 2762.2 & 834.8  \\
    \hline
\multicolumn{5}{l}{}\\
\multicolumn{5}{l}{\text{(c) Observed values for the nontrivial
part of $N_{(\ref{tc})}(100057)$}}\\
        \hline
        g \setminus h & \ANY & \PR & \RP & \RPPR  \\
        \hline
\ANY & 100860 & 9231 & 30291 & 2820\\
\hline
\PR & 30850 & 9231 & 9231 & 2820\\
\hline
\RP &    30368 & 2882 & 9240 & 916\\
\hline
\RPPR & 9376 & 2882 & 2882 & 916\\

    \hline
\end{array}$$
\end{table}

\begin{table}
    \caption{Relationship between solutions to~(\ref{tc})
and solutions to~(\ref{ha})}
    \label{relationstable}
$$\begin{array}{|c|c|c|c|c|}
    \hline
        a \setminus h & \ANY & \PR & \RP & \RPPR  \\
        \hline
    \ANY & & & g \ANY & g \PR \\
     & & &  h \RP &  h \RPPR \\
        \hline
    \PR & & & g \PR & g \PR\\
        & & &  h \RP &  h \RPPR \\
    \hline
    \RP & h \ANY & g \PR & h \RP & g \PR \\
    &  g \ORD h &  h \PR &  g \ORD h & h \RPPR \\

    \hline
    \RPPR & g \PR & g \PR & g \PR & g \PR\\
    &  h \RPPR &  h \RPPR &  h \RPPR &  h \RPPR\\
        \hline

    \end{array}$$
\end{table}


    Finally, the analysis in the last paragraph of
    Section~\ref{sectiontc} is incomplete.  If $\mbox{$\gcd(h,a,p-1)$} = 1$,
    then the correspondence between solutions of~(\ref{tc}) and
    solutions of~(\ref{ha}) is one-to-one.  (E.g., if $h \RP$ or $a
    \RP$.)  If $\gcd(h,a,p-1) > 1$, however, more than one solution
    to~(\ref{tc}) may give the same solution to~(\ref{ha}).  This
    will be explored in more detail in a forthcoming paper.


\end{document}